%% file: epcen.tex
\documentclass[12pt,twoside,a4paper,reqno]{amsart}
\usepackage{amsmath}
\usepackage{amssymb}
\usepackage{graphicx}

\usepackage{graphicx}
\graphicspath{{figs/}}
\usepackage{color}
\numberwithin{figure}{section}

\def\sideremark#1{\ifvmode\leavevmode\fi\vadjust{\vbox to0pt{\vss 
      \hbox to 0pt{\hskip\hsize\hskip1em           
 \vbox{\hsize3cm\tiny\raggedright\pretolerance10000
 \noindent #1\hfill}\hss}\vbox to8pt{\vfil}\vss}}}%

\newtheorem{theorem}{Theorem}
\newtheorem{corollary}{Corollary}
\newtheorem{proposition}{Proposition}
\newtheorem{lemma}{Lemma}

\theoremstyle{definition}

\newtheorem{example}{Example}
\theoremstyle{remark}
\newtheorem{remark}{Remark}

\renewcommand{\epsilon}{\varepsilon}

\newcommand{\var}{\mathcal{V}ar}

\newcommand{\R}{\mathbb{R}}
\newcommand{\C}{\mathbb{C}}

\newcommand{\dr}{\mathrm{d}}
\newcommand{\bbC}{\mathbb{C}}
\newcommand{\bbZ}{\mathbb{Z}}

\newcommand{\bbr}{\mathbb{R}}

\newcommand{\ep}{\epsilon}

\title[pseudo-Abelian integrals: slow-fast case]{Pseudo-Abelian integrals on slow-fast Darboux systems}
\author{Marcin Bobie\'nski}
\address{Institute of Mathematics, Warsaw University,
ul. Banacha 2, 02-097 Warsaw, Poland}
\email{mbobi@mimuw.edu.pl}
\thanks{This research was supported by ANR ANAR Project, Polish MNiSzW Grant No N N201 397937, and by Israel Science Foundation (grant No. 1501/08).}

\author{Pavao Marde\v{s}i\'c}
\address{Universit\'e de Bourgogne, Institut de
Math\'ematiques de Bourgogne, U.M.R. 5584 du C.N.R.S., B.P. 47870, 21078
DIJON CEDEX - FRANCE}
\email{mardesic@u-bourgogne.fr}

\author{Dmitry Novikov}
\address{Department of Mathematics, Weizmann Institute of Science,
Rehovot, Israel}\email{dmitry.novikov@weizmann.ac.il}
\begin{document}
\begin{abstract}
We study pseudo-Abelian integrals associated with polynomial deformations of slow-fast Darboux integrable systems.
Under some assumptions we prove local boundedness of the number of their zeros.
\end{abstract}

\today
\maketitle

\section{Introduction and main result}
 Pseudo Abelian integrals appear as the principal (linear) part of the displacement function in polynomial deformations of Darboux integrable cases. This paper is a part of the program of proving uniform finiteness of the number of zeros of pseudo Abelian integrals, see  \cite{n,bm,bmn}. After studying the generic cases \cite{n,bm}, nongeneric cases must be studied, too. Here we study zeros of pseudo Abelian integrals associated to deformations of slow-fast Darboux integrable systems.

\begin{figure}[htpb]
\input{figs/singfree.pspdftex}
\caption{Phase portrait of the integrable system.}
\label{fig:hatTor}
\end{figure}
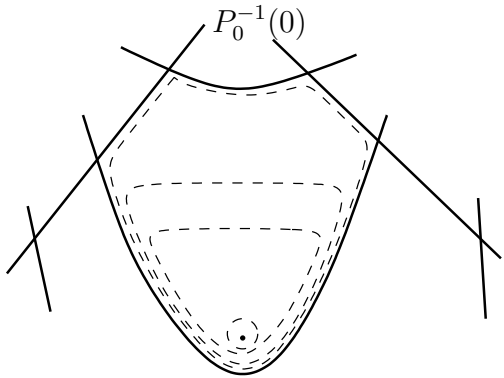

More precisely consider Darboux integrable system $\omega_0$ given by
\begin{equation}\label{omega0}
\omega_0=M\frac{dH_0}{H_0},
\end{equation}
where
\begin{equation}\label{H0}
M=\prod_{i=1}^k P_i, \quad H_0=\prod_{i=1}^k P_i^{a_i},\quad a_i>0,\quad  P_i\in\R[x,y].
\end{equation}
 We consider the family of forms  $\omega_\ep$ given by

\begin{equation}
\label{fol}
\omega_\ep = P_0 M \frac{\dr H_0}{H_0} + \ep M \dr P_0,\quad  P_0\in\R[x,y].
\end{equation}
Note that for $\epsilon=0$ this form defines a foliation singular along  the curve $P_0=0$.
The form \eqref{fol} is Darboux integrable with first integral
\begin{equation}
\label{dbxfn}
H_\ep =  H_0\, P_0^\ep.
\end{equation}
The system \eqref{fol} is a slow-fast and $P_0=0$ is the slow manifold. The fast dynamics is given by the Darboux integrable system \eqref{omega0}.
Figure \ref{fig:hatTor} presents typical example that we treat.

Assume that the system \eqref{fol} has a family $\gamma_{\ep}(t)\subset H_\ep^{-1}(t)$ of cycles.
Consider the polynomial perturbation of the system \eqref{fol} .
\begin{equation}
	\label{torontopert}
	\omega_{\ep,\delta} = \omega_\ep + \epsilon_1 \eta, \quad\epsilon_1>0.
\end{equation}

The linearization in perturbation parameter $\epsilon_1$ of the Poincar\'e first return map is given by the pseudo-Abelian integral
\begin{equation}
	\label{torint}
	I_\ep(t) = \int_{\gamma_\ep(t)} \frac{\eta}{P_0 M}.
\end{equation}

In this paper we study pseudo-Abelian integrals taken along the cycles $\gamma_\ep(t)$ of the simplest
integrable system bifurcating from a slow-fast system:

Assume that a compact region $D$ is  bounded by $P_0=0$ and some separatrices $P_i=0$, $i=1,\ldots,k$.
Assume that the functions $P_i$, $i=0,\ldots k$  are smooth and intersect transversally
in $D$ and  that the foliation \eqref{omega0} has no singularities on $\mathrm{Int} D$.

Assume moreover that  $P_0=0$ is transverse to the foliation $\omega_0=0$ in all points of $D\cap \{P_0=0\}$ except for one point $p_0$, where the contact is quadratic.
Then for $\ep\ne0$ a singular point $p_\ep$ bifurcates from $p_0$.
It corresponds to a real value $t_{\ep}=H_\ep(p_\ep)$.
The bifurcating singular  point $p_{\ep}$ is a center entering the domain $D$ for $\ep>0$.
Let $\gamma_\ep(t)$ be the family of cycles in the basin of the center bifurcating from $p_0$. For each $\ep$, the cycles are defined on an interval
$[0,t_{{\ep}}]$.
\begin{remark}
Note that if $D$ is compact, then necessarily  $\{P_0=0\}$ must have at least one contact point with the foliation $\omega_0=0$.
Below in the paper we call \emph{turning point} the point of tangency between $\{P_0=0\}$ and the leaves of $\omega_0=0$.

For $\ep<0$ instead of a center, a saddle point bifurcates from $p_0$ just outside of $D$.
\end{remark}
\begin{example}\label{ex:ttoronto}
Consider  the foliation
\begin{equation}	\label{toronto}
	\omega_\ep = P_0 \dr P_1 + \ep P_1 \dr P_0,\qquad P_0=y-x^2,\quad P_1=1-y,
\end{equation}
with the first integral $H_\ep=(1-y)(y-x^2)^{\epsilon}$. It has a critical point $p_\ep=(\frac {\ep}{1+\ep},0)$,  which is a saddle for $\epsilon<0$, coincides with the tangency point $p_0=(0,0)$ for $\ep=0$, and is a center for $\ep>0$. For each $\ep>0$ the region $D=\{x^2\le y\le 1\}\subset \R^2$ is filled by a nest of cycles $\gamma_\ep(t)$ vanishing at $p_\ep$. These cycles are parameterized by values of $t$ varying from $0$ (corresponds to the polycycle $\partial D$) to 
$$
t_\ep=H_\ep(p_\ep)=(1+\ep)^{-1}\left(1+\frac 1 \ep\right)^{-\ep}\to e^{-1} \quad\text{as}\quad \ep\to 0.
$$

\begin{figure}[htpb]
\includegraphics[width=0.9\hsize]{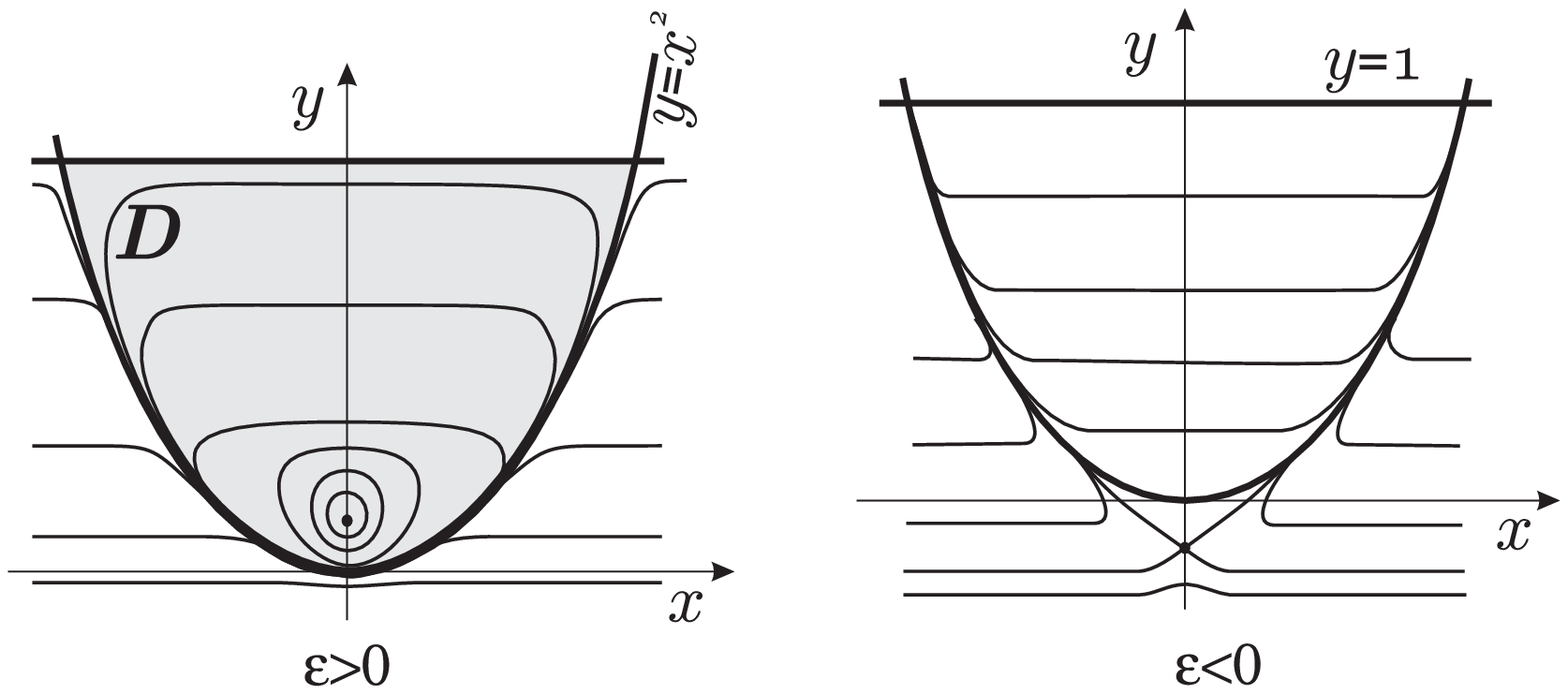}
\caption{Phase portrait of $H_\ep=(1-y)(y-x^2)^{\ep}$}
\label{fig:TTor}
\end{figure}

\end{example}

\begin{theorem}\label{main}
Let $I_\ep(t)$ be the family of pseudo-Abelian integrals as defined above with the above genericity assumptions. Then
there exists a bound  for the number of isolated zeros of the pseudo-Abelian integrals $I(t,\ep)$, for $\ep>O$ and $t\in(0,t_{\ep})$. The bound is locally uniform with respect to all parameters: the form varying in some finite dimensional analytic family, the coefficients of the polynomials, exponents  and in particular with respect to $\ep$.
\end{theorem}

\begin{remark}
The estimates on the number of zeros of pseudo-Abelian integrals do not imply in general estimates on cyclicity of corresponding polycycles, due to presence of the so-called "alien" cycles, see e.g. \cite{DR2}.
\end{remark}

In order to prove the theorem we want to apply variation and argument principle as in \cite{bm,b,bmn}.
The difficulty lies in the fact that  the family of cycles of integration $\gamma_\ep(t)$, $\ep>0$ does not have a smooth limit as $\ep\to0.$ The family tends to a family of broken cycles formed by  parts of the slow manifold and parts of leaves of the fast foliation.  After taking an $\ep$-scaled variation,
the cycle of integration $\gamma_\ep(t)$ is replaced by a \emph{figure eight cycle} $\delta_\ep(t)$, see Proposition~\ref{pr:vartor0}, and the integral $I_\ep(t)$ by the integral $J_\ep(t)$ taken along $\delta_\ep(t)$.
For $t$ not too close to the center value $t_\ep$, the cycle $\delta_\ep(t)$ can be moved along leaves of the foliation to be at some positive $\ep$-independent distance from the slow manifold. This property, which does not hold for the initial cycle $\gamma_\ep(t)$, implies good analytic properties of $J$.

However, as $t\to t_\ep$, and $\ep\to0$, both cycles $\gamma_\ep(t)$ and $\delta_\ep(t)$ approach the turning point. To overcome  this difficulty in a $\ep$-uniform way, we perform a blow-up in the family $\omega_\ep$. To each chart  of the blow up corresponds a time scale. In blown-up coordinates the cycles $\gamma_\ep(t)$ and $\delta_\ep(t)$ have smooth limits in respective charts. This proves analyticity of $I$ or  $J$ in the timescale of the convenient chart.

\subsection{Acknowledgment}
The authors express their gratitude to Institute of Mathematics of Warsaw University and Weizmann Institute of Science for hospitality.

\section{$\ep$-variation of the $\gamma_\ep(t)$ cycle}
\label{sec:toronto}

The key point of our approach is to study analytic properties of the integral along the family of \emph{figure eight loops}.

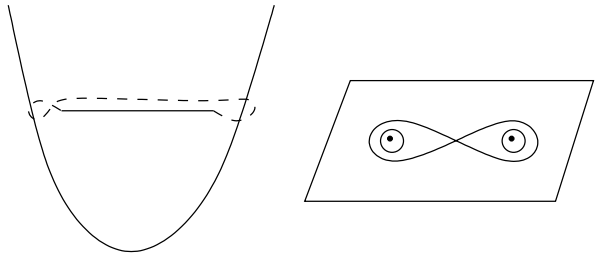
\begin{figure}[htpb]
\input{figs/tor8.pspdftex}
\caption{The "figure 8" loop}
\label{fig:tor8}
\end{figure}

On the smooth leaf $H_0=t$ of the unperturbed equation there exists a figure eight cycle lying close to the real segment $\{H_0=t\}\cap D$ and winding around the points of intersection of this leaf with the curve $P_0=0$, see Figure \ref{fig:tor8}. Since this cycle lies on a finite distance from this curve, i.e. in a domain where the leaves of the foliation $H_\epsilon=t$ depend analytically on $\epsilon, t$, one can move it to a complex leaf of the foliation $H_\ep=t$, uniquely up to a small homotopy.
The figure eight family of cycles $\delta_\ep(t)$ is  defined as the family of these lifts.

Let
 $\gamma_{\ep}(t)\subset H_\ep^{-1}(t)$
as above and $I_\ep(t)$ be the integral \eqref{torint} along $\gamma_\ep(t)$.

\begin{proposition}\label{pr:vartor0}
The $\ep$-variation of the cycle $\gamma_\ep(t)$  is equal to the figure eight cycle $\delta_\ep(t)$:
\begin{equation}
\var_\ep(\gamma_\ep(t))=\gamma_\ep(t\ e^{\pi i \ep})-\gamma_\epsilon(t\ e^{-\pi i \ep})=\delta_\epsilon(t).
\end{equation}
\end{proposition}

\begin{proof} As in \cite{bm, bmn}, we transport the cycle $\gamma_\ep(t)$ for $t\mapsto t\ e^{\pi i \ep}$ and subtract the transport along $t\mapsto t\ e^{-\pi i \ep}. $ Using the local expression of $H_\ep$ near the slow manifold, we see that the parts of  $\gamma_\ep(t\ e^{\pi i \ep})$ and $\gamma_\ep(t\ e^{-\pi i \ep})$ along the slow manifold cancel.
The figure eight cycle $\delta_\ep(t)$ remains.
\end{proof}

\begin{corollary}
\label{pr:vartor}
The $\ep$-variation of the pseudo-Abelian integral $I_\ep(t)$ is an integral of the form $\frac{\eta}{P_0 M}$ along the figure eight cycle $\delta_\ep(t)$.
\begin{equation}
	\label{vartor}
(\var_\ep I) (t) \colon= I_\ep(t\ e^{\pi i \ep}) - I_\ep(t\ e^{-\pi i \ep}) = \int_{\delta_\ep(t)}\frac{\eta}{P_0 M} = J_\ep (t),
\end{equation}

\end{corollary}

\section{Blowing up the turning point}\label{sssec:blowup}
We need to prove analytic properties of the integrals $I$ and $J$ in a neighborhood of the slow manifold with respect to a convenient time scale. This will be achieved by a convenient blow-up of the family in a neighborhood of the turning point. This section is dedicated to this geometric construction.

Note that by Morse lemma we can put our family \eqref{fol} to a normal form of Example~\ref{ex:ttoronto} in a neighborhood of the turning point:

\begin{lemma}\label{normal}
By an analytic $\ep$-independent change of coordinates defined in a neighborhood of the turning point $p_0$ we can assume that $\omega_0=dy$ and $P_0=y-x^2$.
\end{lemma}

We want to study the analytical properties of the foliation in a  neighborhood of the turning point (the vertex of parabola).
The difficulty is the approaching of the center to the vertex of the parabola, linearly with $\ep$. So we make a blow-up of the of the turning point of the family in the product space $(x,y,\ep)$ of phase and parameter spaces.  The family blow-ups were studied by Dumortier and Roussarie e.g. in \cite{RRRR}. This is needed since we want to prove analyticity with respect to both phase and parameter values. We want our blow-up to preserve the parabola $y=x^2$ and to separate the newborn center from the vertex. This requirements lead to the weighted blow-up with weights $(1:2:2)$.

Recall the construction of this weighted blow-up. We define the weighted projective space $\C P^2_{1:2:2}$ as the factor space of $\C^3$ by the $\C^*$ action $(x,y,\ep)\mapsto(tx, t^2y, t^2\ep)$. The blow-up of $\C^3$ at the origin is defined as the incidence three dimensional manifold $M=\{(p,q)| q\in p\}\subset \C P^2_{1:2:2}\times \C^3$.

Note that the weighted action of $\bbC^*$ on $\bbC^3\setminus \{0\}$ is not free. The stabilizer of points $(0,y,\ep)$ is a $\bbZ/2=\{\pm 1\}$ subgroup. As a result, the quotient is not smooth on the line  $\{x=0\}$. In chart described below we will work with the double covering which is smooth.

The blow-down $\pi:M\to\C^3$ is just the restriction to $M$ of the projection $\C P^2_{1:2:2}\times \C^3\to\C^3$.

For future applications we will need explicit formulae  for the blow-up in the standard affine charts of $M$.
The projective space $\C P^2_{1:2:2}$ is covered by three affine charts: $U_1=\{x\not=0\}$ with coordinates $(Y_1,E_1)$,  $U_2=\{y\not=0\}$ with coordinates $(X_2,E_2)$ and  $U_3=\{\ep\not=0\}$ with coordinates $(X_3,Y_3)$.  The transition formulae follow from the requirement that the points
$(1,Y_1,E_1)$, $(X_2,1,E_2)$ and $(X_3,Y_3,1)$ lie on the same orbit of the action:
\begin{equation}
\begin{split}
(Y_1,E_1)\mapsto (X_2=\frac 1 {\sqrt{Y_1}}, E_2=\frac{E_1}{Y_1})\\
(Y_1,E_1)\mapsto (X_3=\frac 1 {\sqrt{E_1}}, Y_3=\frac{Y_1}{E_1}).
\end{split}
\end{equation}

\begin{remark}
Note that transition functions are singular since we cover smooth double covering instead of the exceptional divisor itself. It is easy to observe that replacing coordinates $X_2$ and $X_3$ in charts $U_2$ and $U_3$ respectively, we obtain usual, smooth projective transition functions.
\end{remark}

These affine charts define affine charts on $M$, with coordinates $(Y_1,E_1,t_1)$, $(X_2,E_2, t_2)$ and $(X_3,Y_3,t_3)$. The projection $\pi$ is written as
\begin{align}
&x=t_1  & y=t_1^2 Y_1 \quad& \ep=t_1^2 E_1&\\
&x=t_2X_2 & y=t_2^2 \quad& \ep=t_2^2 E_2&\\
&x=t_3X_3 & y=t_3^2Y_3 \quad& \ep=t_3^2.&
\end{align}

We apply this blow-up $\pi$  to the one-dimensional foliation $\mathcal{F}$  on $\C^3$ given by the intersection of $d\ep=0$ and $\omega_{\ep}=0$.
This foliation has a complicated singularity at the origin. Denote by $\pi^{-1}\mathcal{F}$ the lifting of the foliation $\mathcal{F}$ to the complement of the exceptional divisor $\pi^{-1}(0)$. This foliation is regular outside of the preimage of the parabola $\mathcal{P}=\{y=x^2, \ep=0\}$.
\begin{proposition}\label{pr:blowup}
The foliation $\pi^{-1}\mathcal{F}$  can be extended analytically to the exceptional divisor $\pi^{-1}(0)$. The resulting foliation $\pi^*\mathcal{F}$ is regular outside of the strict transform of the parabola $\mathcal{P}$, the family of centers $\{X_3=0, Y_3=\frac{1}{1+t_3^2}\}$ and a point $p_s=(1:0:0)$ on the exceptional divisor.
\end{proposition}
\begin{figure}[htpb]
\includegraphics[width=0.5\hsize]{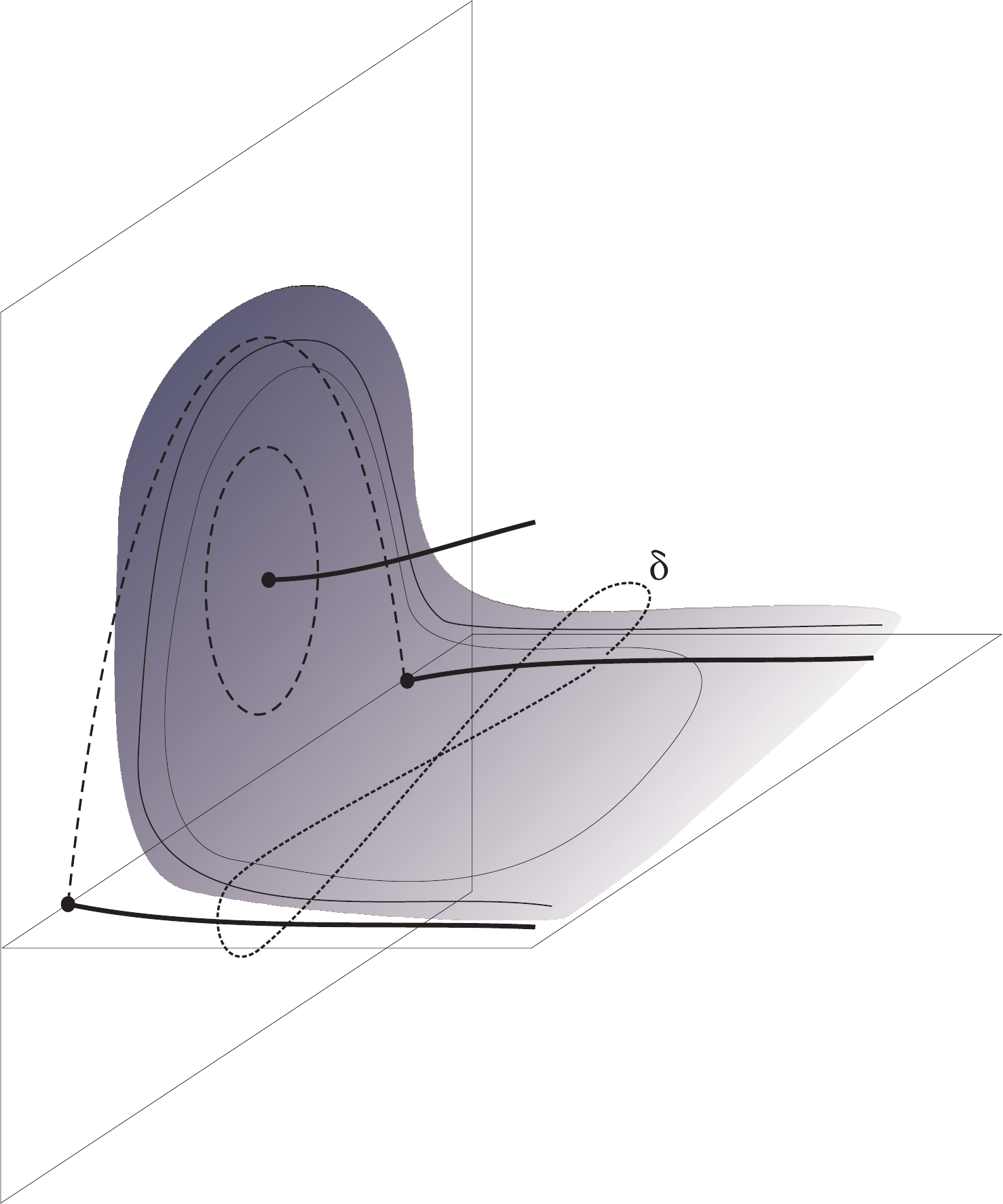}
\caption{The foliation $\pi^{-1}\mathcal{F}$}
\label{fig:blup}
\end{figure}

\begin{remark}
The additional singular point has clear geometric interpretation which is characteristic to the family blow-up of a slow-fast system. It is a "trace" of the fast direction on the blow-up of the slow manifold. Consider the following toy example of foliation in $\bbr^2$ defined by 1-form $\omega_\ep = x \dr (1-y+ax) + \ep (1-y+ax)\dr x$. Making the family blow up
\[
\pi:\qquad y=x\, Y,\quad \ep = x\, E
\]
we obtain a foliation given by 2-form
\[
\tfrac{1}{x^2}\, \pi^*(\omega_\ep\wedge \dr \ep) = (-a -E -a E x + Y + E x Y)\dr E\wedge \dr x + x\,\dr E \wedge \dr Y+ E \dr x\wedge \dr Y
\]
that vanishes at the point $x=0, E=0, Y=a$ corresponding to the direction of the fast system.
\end{remark}

\begin{proof} We check it in each chart separately. Note that a codimension 2 foliation $\mathcal{F}$ in 3-dimensional space is uniquely defined by 2-form $\sigma$ by the condition
\[
\mathcal{F}=\{v:\ i_v \sigma =0\}.
\]
Alternatively, locally $\sigma=i_v \operatorname{Vol}$, where $\operatorname{Vol}$ is a suitable non-vanishing volume form.
Singular points of foliation correspond to zeros of 2-form $\sigma$. The foliation \eqref{toronto} in 3-dimensional space $\bbC^3$ is given by $\sigma_\ep = \omega_\ep\wedge \dr \ep$. The pull-back foliation $\pi^* \mathcal{F}$ (strict transform of $\mathcal{F}$) is defined by the pull-back $\pi^* \sigma_\ep$ divided by a suitable power of the function defining the exceptional divisor. In charts $U_j,\ j=1,2,3$ we have $\pi^* \sigma = t_j^5 \widetilde{\sigma}_j$, where
\begin{align}\label{eq:blownup fol}
\widetilde{\sigma}_1 =&
2 (Y_1-1)(Y_1-E_1+E_1 t_1^2 Y_1)\, \dr E_1\wedge \dr t_1 + \nonumber\\
&+(Y_1-1-E_1+E_1 t_1^2 Y_1) \Big(t_1 \dr E_1\wedge \dr Y_1 + 2 E_1\dr t_1\wedge \dr Y_1\Big),\nonumber\\
\widetilde{\sigma}_2 =& 2 (1-X_2^2)(1-E_2+E_2 t_2^2)\,\dr E_2\wedge\dr t_2+ \\
&+ 2(1-t_2^2) X_2E_2\Big( t_2 \dr E_2\wedge \dr X_2+ 2 E_2 \dr t_2\wedge \dr X_2\Big),\nonumber\\
\widetilde{\sigma}_3 =& 4 X_3 (1-t_3^2 Y_3)\dr t_3\wedge \dr X_3 - 2 (1+X_3^2-Y_3-t_3^2 Y_3)\dr t_3\wedge \dr Y_3\nonumber
\end{align}
The zero locus of the form $\widetilde{\sigma}$ in a neighborhood of the exceptional divisor consists of germs of two curves and the singular point $p_s=[1:0:0]$ generated by weighted action. These curves are $(X_2=\pm 1,E_2=0)$ (strict transform of the parabola $P=0$) and $(X_3=0,Y_3=\tfrac{1}{1+t_3^2})$ (family of centers).

In chart $U_3$ we observe that the lifted foliation is given by two first integrals: $t_3$ and
\[
s=(1-t_3^2 Y_3)^{1/t_3^2} (Y_3-X_3^2)
\]
which can be analytically continued to $t_3= 0$. This foliation has no singularities near the exceptional divisor except for the line of centers $X_3=0, Y_3=\frac{1}{1+t_3^2}$. Note that the strict transform of the parabola is outside of this chart.

\end{proof}

\section{Proof of the Theorem}

In this section we first take benefit from the blowing-up in the family performed in the previous section to prove analyticity of the integrals $I$ and $J$ in convenient time scale. Let us be given a compact family of cycles on the exceptional divisor of the blown-up foliation at  a finite distance from the singularities.
We can extend it to a continuous family in its full neighborhood in the total blown-up space. By analyticity of blown-up foliation, see Proposition~\ref{pr:blowup}, the integrals along these cycles will depend analytically on the cycle in the extended family.
In particular, this means that:
\begin{lemma}\label{anal} For any $R_s < s_c(0)=e^{-1}$
\begin{enumerate}
\item  the integral $I$ is an analytic  function of $s$ and $t_3$ in some neighborhood of the arc  $\{0\}\times\{R_s e^{i\varphi}||\varphi|\le \pi\}$ in complex $(t_3, s)$-plane;
\item the integral $I$ is an analytic function of $s$ and $t_3$ in some neighborhood of the segment $\{0\}\times [R_s, s_c(0)]$  in complex $(t_3, s)$-plane;
\end{enumerate}
\end{lemma}
\begin{proof}[Proof of Lemma \ref{anal}]
Note that the restriction of the foliation $\pi^*\mathcal{F}$ to the exceptional divisor has first integral $s=e^{-Y_3}\left( Y_3-X_3^2\right)$. It is easy to construct a compact family of cycles $\gamma$ lying on the exceptional divisor and corresponding to the values of $s$ mentioned in the first claim of the Lemma, so above argument proves the first claim.

For the second claim consider the compact family of real cycles vanishing at the center. Analyticity of the integral along these cycles is standard.
\end{proof}

Consider the transversal $X_2=0$ to the line $t_2=E_2=0$ in chart $U_2$. To each point on this transversal corresponds a figure eight cycle as in Proposition~\ref{pr:vartor0} passing through this point and lying on a leaf of the foliation $\pi^*\mathcal{F}$ (note that this line is a leaf of $\pi^*\mathcal{F}$), so the integral $J$ defines  a function on this transversal in a neighborhood of $0$.
\begin{lemma}\label{anal2}
 The function $J$  is an analytic function of  $t_2, E_2$ for $|E_2|, |t_2|\le \rho$ for some sufficiently small $\rho>0$.
\end{lemma}

\begin{proof}
On the line $t_2=E_2=0$ lying in $U_2$ the blown-up foliation \eqref{eq:blownup fol} has two singular points $X_2=\pm1$, see Figure~\ref{fig:blup}. The figure eight cycle lying on this line is on finite distance from these two singularities. By integrability of $\pi^*\mathcal{F}$,  it can be extended to all sufficiently close leaves, forming a continuous family of figure eight cycles.

\end{proof}

\begin{proof}[Proof of Theorem \ref{main}]
The first integral maps the open nest of vanishing cycles to the interval $(0,t_\ep)$. We split the interval into two parts, one being the image $l_0$ of the interval $[R_s, s_c)$, and the remaining part. On the first part the number of zeros of $I$ is uniformly bounded by Gabrielov's theorem, due to Lemma~\ref{anal}.(2).

To estimate the number of zeros on the remaining part we apply argument principle to the contour $\Gamma$ consisting of two arcs $A_\ep$  and $a_\ep$ and two straight segments $L_{\pm}$ joining their ends, see Figure~\ref{fig:argtor}. The arc $a_\ep$ of angle $2\ep$ has infinitesimally  small radius and the arc $A_\ep$ is the image of the arc described in Lemma~\ref{anal}.(1).

The $\ep$-uniform bound for the increment of argument along $A_\ep$ is a direct consequence of the Lemma \ref{anal}.(2).

The increment of argument along segments $L_{\pm}$ can be estimated by number of zeros of the imaginary part of the function $I$ on $L_{\pm}$. We calculate using fact that function $I$ is real on the real segment:
\begin{multline*}
\frac 1 \pi\Delta Arg_{L_+} I \leq \ \#\{t: Im I =0\} +1 = \ \#\{t: I - \overline {I} =0\} +1 =\\
 \ \#\{t: I(e^{\pi i \ep}\, t)- I(e^{-\pi i \ep}\, t)  =0\}+1 =  \, \#\{t: J (t) =0\}+1.
\end{multline*}
Thus, one translates problem into estimating the number of zeros of the integral $J$ along the family of  figure eight cycles corresponding to segments $L_\pm$.

We split the segment $L_\pm$ into two parts. One, closer to $t_c$, is in the image of the transversal considered in Lemma~\ref{anal2}. The functions $J$ on this part can be considered as an analytic function of $t_2, E_2$. The remaining part of $L_\pm$  corresponds to the part of the family of figure eight cycles lying not closer than some fixed positive distance from the vertex of the parabola. These cycles can be deformed along leaves to be not closer than some positive constant from the parabola, see Figure~\ref{fig:tor8}. Therefore the function $J$ is an analytic in $t,\ep$ on this part.

By Gabrielov theorem, the number of zeros of $J$ on $L_\pm$ is uniformly bounded in $\ep$.

To obtain an upper bound for increment of argument of $I$  along the small arc $a_\ep$ we use an argument similar  to those from \cite{bm,bmn,b,n}. One can easily prove that $|I|\le C|t|^{M/\ep}$ in sectors.  On the other hand, for any fixed $\ep$ the function $I$ has a leading term at $t=0$ of the form $t^\beta \log ^j t$, see the above papers. Together it proves the existence of the uniform in $\ep$ upper bound for the increment of argument.

All the above constructions depend analytically on parameters like coefficients of the polynomials $P_i$, exponents $a_i$ and coefficients of the form $\eta$.

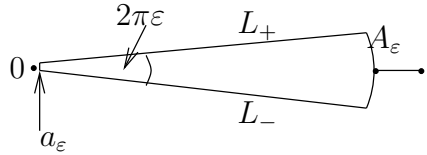
\begin{figure}[htpb]
\input{figs/argtor2.pspdftex}
\caption{Contour $\Gamma$}
\label{fig:argtor}
\end{figure}

\end{proof}

\section{Concluding remarks and open problems}

\begin{remark}
Consider the  model Example~\ref{ex:ttoronto}.
An alternative way to prove Theorem~\ref{main} for this system would be to reduce it to the situation considered in our previous paper \cite{bmn}.
Indeed, change of variables $x=\sqrt{\ep}X, y=\ep Y$ the system \eqref{toronto} has the first integral $s=(Y-X^2)(1-\ep Y)^{1/\ep}$. After additional blow-up of the point at infinity, we get an unfolding of a system with two saddle-nodes. Unfolding of systems with polycycle with two saddle-nodes was considered in \cite{bmn}.
However, this approach does not generalizes to more general systems considered in Theorem~\ref{main}.
\end{remark}

In the spirit of our program of proving of uniform finiteness of the number of zeros of pseudo-Abelian integrals, more general slow-fast Darboux systems should be studied. Already in the generic, stable under small perturbations of coefficients  cases one has to face several problems, in particular:
\begin{enumerate}
\item The unperturbed Darboux system $\omega_0=0$ can have extra singular points inside $D$, not lying on the zero level of $H_0$;
\item Unperturbed Darboux system can have a nest of cycles, and the slow manifold  cuts the nest and becomes a part of polycycle bounding the new nest;
\item The curve $\{P_0=0\}$ can have additional tangency points with the leaves of $\omega_0=0$, which generate saddles type singularities;
\item The nest of cycles accumulating to the polycycle and encircling more than one newborn singular point, so called "big cycles";
\end{enumerate}

These possible scenarios are illustrated on Figure~\ref{fig:biftor} below, keeping numeration.

\begin{figure}[htpb]
\input{figs/dbxscasescal.pspdftex} \input{figs/ph2scaled.pspdftex}
\caption{Phase portrait of the system with two factors.}
\label{fig:biftor}
\end{figure}
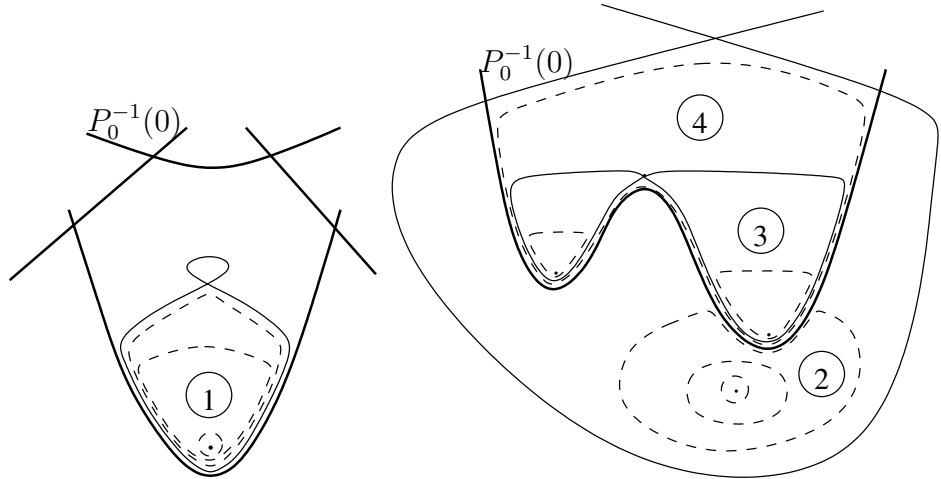

\end{document}

%% file: figs/singfree.pspdftex
\begin{picture}(0,0)%
\includegraphics{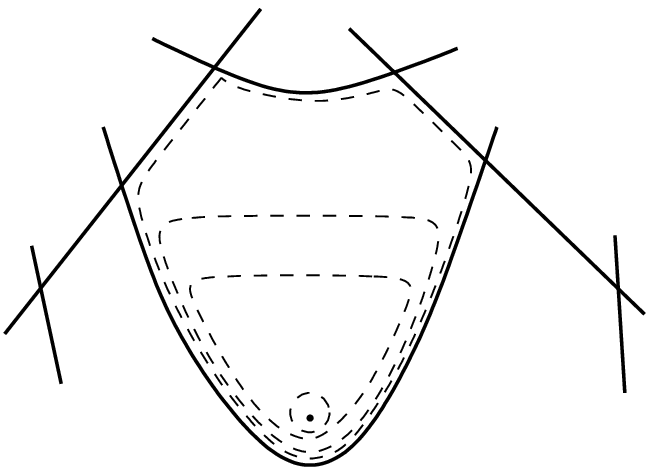}%
\end{picture}%
\setlength{\unitlength}{4144sp}%
\begingroup\makeatletter\ifx\SetFigFontNFSS\undefined%
\gdef\SetFigFontNFSS#1#2#3#4#5{%
  \reset@font\fontsize{#1}{#2pt}%
  \fontfamily{#3}\fontseries{#4}\fontshape{#5}%
  \selectfont}%
\fi\endgroup%
\begin{picture}(2969,2221)(6324,-3023)
\put(7561,-961){\makebox(0,0)[lb]{\smash{{\SetFigFontNFSS{12}{14.4}{\rmdefault}{\mddefault}{\updefault}{\color[rgb]{0,0,0}$P_0^{-1}(0)$}%
}}}}
\end{picture}%

%% file: figs/tor8.pspdftex
\begin{picture}(0,0)%
\includegraphics{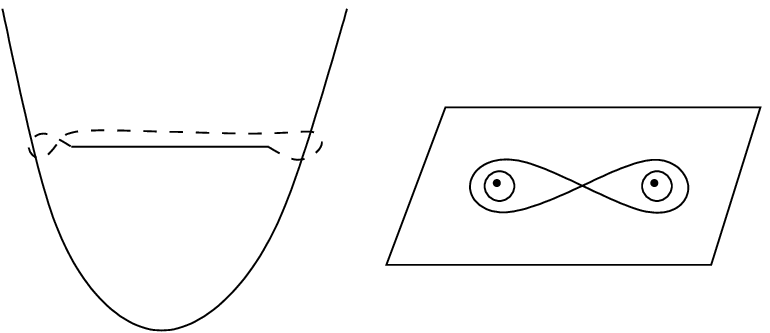}%
\end{picture}%
\setlength{\unitlength}{4144sp}%
\begingroup\makeatletter\ifx\SetFigFontNFSS\undefined%
\gdef\SetFigFontNFSS#1#2#3#4#5{%
  \reset@font\fontsize{#1}{#2pt}%
  \fontfamily{#3}\fontseries{#4}\fontshape{#5}%
  \selectfont}%
\fi\endgroup%
\begin{picture}(3489,1494)(34,-1903)
\end{picture}%

%% file: figs/argtor2.pspdftex
\begin{picture}(0,0)%
\includegraphics{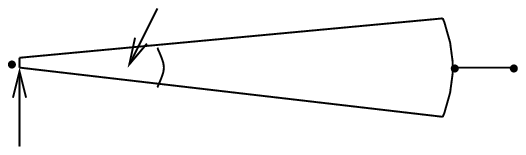}%
\end{picture}%
\setlength{\unitlength}{4144sp}%
\begingroup\makeatletter\ifx\SetFigFontNFSS\undefined%
\gdef\SetFigFontNFSS#1#2#3#4#5{%
  \reset@font\fontsize{#1}{#2pt}%
  \fontfamily{#3}\fontseries{#4}\fontshape{#5}%
  \selectfont}%
\fi\endgroup%
\begin{picture}(2477,964)(526,-1709)
\put(1171,-916){\makebox(0,0)[lb]{\smash{{\SetFigFontNFSS{12}{14.4}{\rmdefault}{\mddefault}{\updefault}{\color[rgb]{0,0,0}$2\pi \ep$}%
}}}}
\put(1891,-961){\makebox(0,0)[lb]{\smash{{\SetFigFontNFSS{12}{14.4}{\rmdefault}{\mddefault}{\updefault}{\color[rgb]{0,0,0}$L_+$}%
}}}}
\put(1891,-1501){\makebox(0,0)[lb]{\smash{{\SetFigFontNFSS{12}{14.4}{\rmdefault}{\mddefault}{\updefault}{\color[rgb]{0,0,0}$L_-$}%
}}}}
\put(721,-1636){\makebox(0,0)[lb]{\smash{{\SetFigFontNFSS{12}{14.4}{\rmdefault}{\mddefault}{\updefault}{\color[rgb]{0,0,0}$a_\ep$}%
}}}}
\put(541,-1231){\makebox(0,0)[lb]{\smash{{\SetFigFontNFSS{12}{14.4}{\rmdefault}{\mddefault}{\updefault}{\color[rgb]{0,0,0}$0$}%
}}}}
\put(2656,-1051){\makebox(0,0)[lb]{\smash{{\SetFigFontNFSS{12}{14.4}{\rmdefault}{\mddefault}{\updefault}{\color[rgb]{0,0,0}$A_\ep$}%
}}}}
\end{picture}%

%% file: figs/dbxscasescal.pspdftex
\begin{picture}(0,0)%
\includegraphics{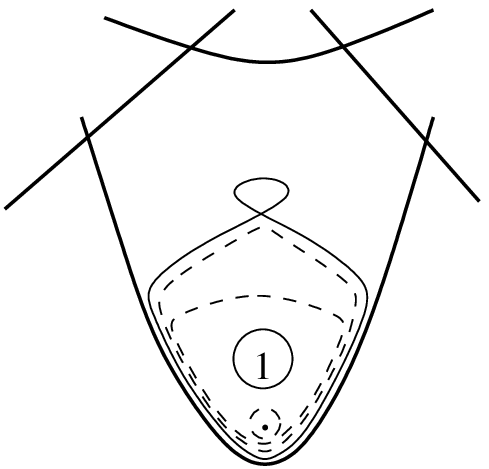}%
\end{picture}%
\setlength{\unitlength}{4144sp}%
\begingroup\makeatletter\ifx\SetFigFontNFSS\undefined%
\gdef\SetFigFontNFSS#1#2#3#4#5{%
  \reset@font\fontsize{#1}{#2pt}%
  \fontfamily{#3}\fontseries{#4}\fontshape{#5}%
  \selectfont}%
\fi\endgroup%
\begin{picture}(2213,2238)(14,-1330)
\put(496,749){\makebox(0,0)[lb]{\smash{{\SetFigFontNFSS{12}{14.4}{\rmdefault}{\mddefault}{\updefault}{\color[rgb]{0,0,0}$P_0^{-1}(0)$}%
}}}}
\end{picture}%

%% file: figs/ph2scaled.pspdftex
\begin{picture}(0,0)%
\includegraphics{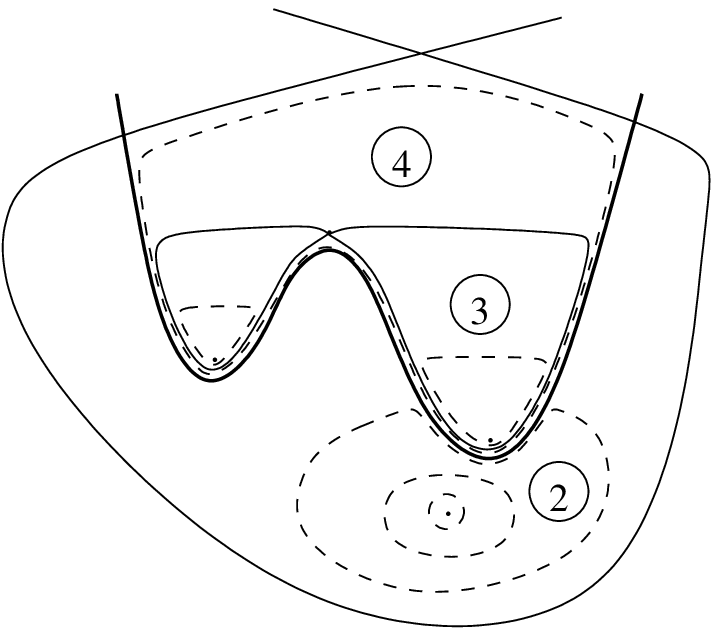}%
\end{picture}%
\setlength{\unitlength}{4144sp}%
\begingroup\makeatletter\ifx\SetFigFontNFSS\undefined%
\gdef\SetFigFontNFSS#1#2#3#4#5{%
  \reset@font\fontsize{#1}{#2pt}%
  \fontfamily{#3}\fontseries{#4}\fontshape{#5}%
  \selectfont}%
\fi\endgroup%
\begin{picture}(3256,2845)(11,-1994)
\put(541,434){\makebox(0,0)[lb]{\smash{{\SetFigFontNFSS{12}{14.4}{\rmdefault}{\mddefault}{\updefault}{\color[rgb]{0,0,0}$P_0^{-1}(0)$}%
}}}}
\end{picture}%